\magnification = 1100

\input amstex
\define\osc{\text{\rm osc}\,}
\documentstyle{amsppt}
\NoBlackBoxes

\topmatter
\title Locally uniformly convex norms in Banach spaces and
their duals
\endtitle
\author Richard Haydon
\endauthor
\affil Brasenose College, Oxford \endaffil
\abstract
It is shown that a Banach space with locally uniformly convex dual
admits an equivalent norm that is itself locally uniformly convex.
\endabstract
\address Brasenose College,
Oxford OX1 4AJ,
U.K
\endaddress
\email richard.haydon\@brasenose.oxford.ac.uk
\endemail
\endtopmatter

\document

\heading
1. Introduction
\endheading

If we consider a real Banach space $Z$ under a norm $\|\cdot\|$ and
its dual space $Z^*$, equipped with the dual norm $\|\cdot\|^*$,
there are important and well-established connections between
convexity properties of $\|\cdot\|^*$ and smoothness properties of
$\|\cdot\|$.  Indeed, strict convexity of $\|\cdot\|^*$ implies
G\^ateaux-smoothness of $\|\cdot\|$, locally uniform convexity of
$\|\cdot\|^*$ implies Fr\'echet-smoothness of $\|\cdot\|$ and
uniform convexity of  $\|\cdot\|^*$ is equivalent to uniform
smoothness of  $\|\cdot\|$.  On the other hand, there would seem to
be, a priori, no reason why a convexity condition in the dual space
$Z^*$ should imply any sort of convexity in $Z$. However, it is a
consequence of the Enflo-Pisier renorming theorem [3,14, or IV.4 of
2] that uniform convexity of $\|\cdot\|^*$ implies that there exists
a norm $|\|\cdot|\|$ on $Z$, {\it equivalent} to the given norm,
which is itself uniformly convex.  One can even arrange that this
new norm  be  both uniformly convex and uniformly smooth.

It is natural to ask whether a similar result about equivalent norms
holds for the weaker properties of strict convexity and locally
uniform convexity.  A counterexample to one of these questions was
given in [8]:  there is a Banach space $Z, \|\cdot\|$ with strictly
convex dual, but such that no equivalent norm on $Z$ is strictly
convex. That the situation may be better for the third property,
locally uniform convexity, was suggested by a theorem of Kenderov
and Moors [10]. This states that a Banach space with locally
uniformly convex dual has the topological property of being
$\sigma$-fragmentable. The main result of the present paper is an
affirmative answer to the full question about locally uniform
convexity.

\proclaim{Theorem A} Let $Z,\|\cdot\|$ be  a Banach space such
that the dual norm $\|\cdot\|^*$ on $Z^*$ is locally uniformly
convex.  There exists an equivalent norm $|\|\cdot|\|$ on $Z$
which is locally uniformly convex. Moreover, $|\|\cdot|\|$ may
be chosen to have  locally uniformly
convex dual norm $|\|\cdot|\|^*$.
\endproclaim

The ``moreover'' statement in Theorem~A is an immediate consequence
of the technique of Asplund averaging, for which the reader is
referred to \S II.4 of [2]. Now it is known [VIII.3.12 of 2] that a
Banach space with a norm which is  locally uniformly convex and has
locally uniformly convex dual norm admits $\Cal C^1$-partitions of
unity: equivalently, on such a space every continuous real-valued
function may be uniformly approximated by functions of class $\Cal
C^1$.  We thus have the following corollary.

\proclaim{Corollary}  Let $Z$ be a Banach space with locally
uniformly convex dual.  Every continuous real-valued function on
$X$ may be uniformly approximated by functions of class $\Cal
C^1$.
\endproclaim

We note that for general Banach spaces $Z$ it is still not known
whether the existence on $Z$ of an equivalent Fr\'echet-smooth norm
(or, more generally, a ``bump function'' of class $\Cal C^1$)
implies $\Cal C^1$ approximability as in the above Corollary.  In
the special case of spaces $Z=\Cal C(K)$, this implication has been
established in [6].

Spaces of the type $\Cal C(K)$ play an important part in our proof
of Theorem~A.  It is of course always the case that we may identify
$Z$ with a subspace of $\Cal C(K)$, where $K$ is the unit ball of
the dual space $Z^*$, equipped with the weak$^*$ topology. When the
dual norm $\|\cdot\|^*$ is locally uniformly convex, this $K$
belongs to what Raja [15] has called the class of {\it
Namioka-Phelps} compacts.  Theorem~A will thus follow from the
following $\Cal C(K)$-renorming theorem.

\proclaim{Theorem B}
Let $K$ be a Namioka-Phelps compact.  Then there is a norm on
$\Cal C(K)$, equivalent to the supremum norm, which is locally
uniformly convex.
\endproclaim

The rest of this paper is devoted to a proof of (a mild
generalization of) Theorem~B.  The definition of a Namioka-Phelps
compact, as well as of the various other topological and renorming
properties with which we are concerned, will be given in the next
section.  We then move on to develop some topological machinery
before  defining a norm in Section~4. The remaining sections contain
the proof that this norm is locally uniformly convex.  The reader
will note the crucial role played by general topology in the proof
that follows: though Theorem~A clearly has some kind of geometrical
content, there is actually surprisingly little geometry in the
proof. The key is the topological concept of a {\it descriptive}
space, due to Hansell [7], and a careful analysis of the
$\sigma$-isolated networks which exist in such spaces. I became
aware of the importance of these notions thanks to the works of the
Hispano-Bulgarian school of geometric functional analysis, and
notably the papers [11, 15,16].  I am particularly grateful to my
friends Anibal Molt\'o and Pepe Orihuela for many helpful
conversations on this material.

\heading
2. Preliminaries
\endheading

Let $Z$ be a real vector space and let $\phi$ be a non-negative
real-valued convex function on $Z$.  When $f\in Z$ and $f_r\in Z$
$(r\in \omega)$, we shall say that the LUR {\it hypothesis holds
for } $\phi$ (and $f$, and the sequence $(f_r)$) if
$$
{\textstyle\frac12}\phi(f)^2+{\textstyle\frac12}\phi(f_r)^2-\phi({\textstyle\frac12}(f+f_r))^2\to
0.
$$
When the function $\phi$ is positively homogeneous, this statement
is equivalent to saying that both $\phi(f_r)$ and
$\phi({\textstyle\frac12}(f+f_r))$ tend to $\phi(f)$ as $r\to
\infty$. This is recorded as Fact II.2.3 in [2], where it is also
noted that, if the function $\phi$ is an $\ell^2$-sum
$\phi^2=\sum_{n=1}^\infty \phi_n^2$ of non-negative convex functions
and if the LUR hypothesis holds for $\phi$, then it holds for each
of the $\phi_n$.  We shall make repeated use of this observation. We
say that a norm $\|\cdot\|$ is  locally uniformly rotund {\it at} a
given element $f$ if, whenever the LUR hypothesis holds for
$\|\cdot\|$, $f$ and a sequence $(f_r)$, we necessarily have
$\|f-f_r\|\to 0$. This brings us back  to a completely standard
definition: we say that a norm on $X$ is {\it locally uniformly
convex} (the term ``locally uniformly rotund'' and its abbreviation
LUR are also used) if it has this property at each $f\in X$.

We now move on to introduce the topological properties that are
relevant to our results. Most of these ideas are due to Hansell [7].
Our terminology follows [16], where a succinct account can be found
of all the results that we need. A crucial notion is that of a {\it
network} for a topology: a collection $\Cal S$ of subsets of $X$ is
said to be a network for the topology $\Cal T$ if every set in $\Cal
T$ is a union of sets in $\Cal S$: that is to say, whenever $x\in
U\in \Cal T$, there exists $N\in \Cal S$ such that $x\in N\subseteq
U$.  A family of sets $\Cal I$ is said to be {\it isolated} for a
topology $\Cal T$ if, for each $N\in \Cal I$, there exists $U\in
\Cal T$ such that $N\subseteq U$ and $U\cap M=\emptyset $ for all
$M\in \Cal I\setminus\{N\}$; equivalently, $N\cap \overline{\bigcup
\Cal I\setminus \{N\}}=\emptyset$.  A family $\Cal S$ is said to be
$\sigma$-isolated if it can be expressed as $\Cal S= \bigcup_{n\in
\omega} \Cal I_n$ with each $\Cal I_n$ isolated.

Let $(X, \Cal T)$ be a topological space and let $d$ be a metric on
$X$ inducing a topology finer than $\Cal T$. We say that the
property $P(d,\Cal T)$ holds if there is a sequence $(B_n)_{n\in
\omega}$ of subsets of $X$ such that the topology generated by $\Cal
T\cup\{B_n:n\in \omega\}$ is finer than the topology $\Cal T_d$
induced by the metric $d$.  An equivalent formulation is that there
exists a sequence $(A_n)_{n\in \omega}$ of subsets of $X$ such that
the intersections $A_n\cap U$, with $U\in \Cal T$, form a network
for $\Cal T_d$. When $P(d,\Cal T)$ holds, there is a network $\Cal
S$ for the metric topology $\Cal T_d$ which is $\sigma$-isolated for
the topology $\Cal T$. An equivalent formulation of this statement
is that, for each $\epsilon >0$, there is a covering $\Cal S$ of
$X$, which is $\sigma$-isolated for $\Cal T$ and which consists of
sets with $d$-diameter at most $\epsilon$. A compact topological
space $(K, \Cal T)$ which has property $P(d, \Cal T)$ for {\it some}
metric $d$ is said to be {\it descriptive}.  There is an intrinsic
characterization of this property:  $K$ is descriptive if and only
if there is a network for $\Cal T$ which is $\Cal
T$-$\sigma$-isolated.  Hansell's general notion of descriptive space
[7] is a space $X$ which is \v Cech-analytic and has a
$\sigma$-isolated network: we are only concerned with descriptive
compact spaces in this paper. Raja [16] shows that the unit ball of
a dual Banach space $Z^*$ is descriptive for its weak* topology if
and only if $Z$ admits an equivalent norm with ``weak* LUR'' dual
norm.

If $(K,\Cal T)$ is compact and has $P(d,\Cal T)$ for some {\it $\Cal
T$-lower semicontinuous} metric $d$, then $K$ is called a {\it
Namioka--Phelps} compact.   Raja [15] has shown that  unit ball of
of a dual Banach space $Z^*$  is a Namioka--Phelps compact (in the
weak* topology) if and only if $Z$ admits an equivalent norm with
LUR dual norm.  The hard part of this theorem is the ``only if''
implication.  In this paper we just use the easy ``if'' implication.

As has already been mentioned in the introduction, we shall obtain
our main theorem from a renorming result for $\Cal C(K)$ where $K$
is a Namioka--Phelps compact.  In fact we prove something  slightly more
general.

\proclaim{Theorem C}
Let $(K, \Cal T)$ be a (descriptive) compact space which has property $P(d,
\Cal T)$ for a  metric $d$.  There is a norm $\|\cdot\|$ on
$\Cal C(K)$, equivalent to the supremum norm, which is locally
uniformly rotund {\it at} $f$, whenever $f$ is both $\Cal
T$-continuous and $d$-uniformly continuous.
\endproclaim

Of course, Theorem C shows that there is a LUR norm on $\Cal C(K)$
provided the metric $d$ can be chosen in such a way that {\it all}
$\Cal T$-continuous functions are $d$-uniformly continuous.  A
metric with this property has been called a Reznichenko metric. It
is easy to see that a lower semi-continuous metric is Reznichenko,
which is why Theorem~C implies Theorem~B.

The reader who is concerned solely with the proofs of theorems A and B
may omit the remainder of this section, which involves
fragmentability, Radon--Nikodym compacta and the delicate
distinction between lower semicontinuous metrics and Reznichenko
metrics.

The topological space $(X, \Cal T)$ is said to be {\it fragmented}
by the metric $d$ if, for every non-empty subset $Y$ of $X$ and
every $\epsilon>0$, there exists $U\in \Cal T$ such that the
intersection $Y\cap U$ is non-empty and of $d$-diameter at most
$\epsilon$.  If there is a sequence $(B_n)_{n\in \omega}$ of subsets
of $X$ such that the topology generated by $\Cal T\cup \{B_n:n\in
\omega\}$ is fragmented by $d$ then we say that $(X,\Cal T)$ is {\it
$\sigma$-fragmented} by $d$.  For more about
$\sigma$-fragmentability the reader is referred to [9]. Clearly
property $P(d,\Cal T)$ implies that $(X,\Cal T)$ is
$\sigma$-fragmented by $d$.  If $(X,\Cal T)$ is descriptive, then
the converse implication holds too.

If $X$ is compact and is fragmented by some lower semicontinuous
metric, we say that $X$ is a {\it Radon-Nikodym compact}. A compact
space is Namioka--Phelps if and only if it is both descriptive and
Radon--Nikodym. The reader is referred to [12] for the basic facts
about this interesting class of spaces.

As has already been remarked, Theorem C leads to a LUR renorming of
$\Cal C(K)$ when $K$ has property $P$ for some Reznichenko metric.
However, it is not clear whether such a compact space also has $P$
for some lower semicontinuous metric.  The situation is closely
related to the open problem of whether every continuous image of a
Radon--Nikodym compact is again a Radon--Nikodym compact. For recent
work on this topic, and the class of {\it quasi-Radon--Nikod\'ym}
compact spaces see [1,5,13].  A compact space is
quasi-Radon--Nikod\'ym if it is fragmented by some Reznichenko
metric.  Every continuous image of a Radon--Nikod\'ym compact is
quasi-Radon--Nikod\'ym.  All this means that we may state a theorem
which may (or may not!) be a generalization of Theorem B as follows.

\proclaim{Theorem D}
If $K$ is descriptive and is a continuous image of a
Radon--Nikodym compact then $\Cal C(K)$ admits a LUR renorming.
\endproclaim

The author does not know whether $\Cal C(K)$ is LUR-renormable for
all descriptive compacta $K$.  By Raja's results the corresponding
question about Banach spaces would be whether a space $Z$ for which
the dual norm on $Z^*$ is w*LUR has itself an equivalent LUR norm.
The most we can get in this direction (using Theorem D, Raja's
theorem and Theorem 1.5.6 of [4]) is the following.

\proclaim{Corollary}
Let $Z$ be a Banach space such that the dual norm on $Z^*$ is
w*LUR.  If, in addition, $Z$ is a subspace of an Asplund-generated space then
$Z$ admits an equivalent LUR norm.
\endproclaim

\heading
3. Descriptive compact spaces and $\sigma$-isolated families
\endheading

The aim of this section is to develop some additional structure in
a descriptive compact space.  We start by making some general observations
about isolated and $\sigma$-isolated families, which are valid without
any compactness assumption.  Let $K$ be a
topological space and let $\Cal I$ be an isolated family of
subsets of $K$.  Then, by definition, we have
$$
N\cap \overline{\bigcup\Cal I\setminus\{N\}} =\emptyset,
$$
for all $N\in \Cal I$.  If we set
$$
\widetilde N =     \overline N\setminus \overline{\bigcup\Cal
I\setminus\{N\}},
$$
and $\widetilde{\Cal I} = \{\widetilde N:N\in \Cal I\}$ then it is clear
that $\widetilde I$ is again an isolated family.  If $N=\widetilde N$ for
all $N\in \Cal I$, we shall say that $\Cal I$ is a {\it regular}
isolated family.

We shall now introduce some notation for regular isolated
families, which will be employed consistently in all that follows.
If $\Cal I$ is a regular isolated family we write $I$ for the
union of the family $\Cal I$, that is
$$
I=\bigcup \Cal I,
$$
and we define
$$
J=\{t\in K: \text{each neighbourhood of } t \text{ meets at least
two members of }\Cal I\}.
$$
By virtue of its definition, $J$ is a closed set.  Moreover, the
closure $\overline I$ is the union of its disjoint subsets $I$ and
$J$; that is to say, $J=\overline I\setminus I$.

Let us now consider a space with a covering $\Cal S$, which is the
union of countably many regular isolated families $\Cal I(i)$
($i\in \omega$).  In accordance with the notation above, we write
$$
I(i)=\bigcup\Cal I(i),\quad J(i)=\overline I(i)\setminus I(i).
$$
We now make a recursive definition of further families $\Cal I(\bold i)=\Cal
I(i_0,\dots,i_{k})$, together with the associated sets $J(\bold
i)$,
when $\bold i= (i_0,\dots, i_{k})\in
\omega^{<\omega}$ is a finite sequence of natural numbers.
$$
\align
I(1_0,\dots,i_k)&=\bigcup\Cal I(1_0,\dots,i_k)\\
J(i_0,\dots,i_{k}) &=  \overline I(i_0,\dots,i_{k})\setminus
I(i_0,\dots,i_{k})\\
\Cal I(i_0,\dots,i_{k},i_{k+1}) &= \{N\cap J(i_0,\dots,i_{k}):
N\in \Cal I(i_{k+1})\}.
\endalign
$$

\proclaim{Lemma 3.1}
If $\bold i =(i_0, \dots, i_k)$ and $0\le l<k$ then
$$
I(i_0,\dots,i_{k})\subseteq J(i_0,\dots,i_{l})\subseteq
J(i_l).
$$
If the natural numbers $i_0,i_1,\dots,i_{k}$ are not all
distinct then $I(i_0,\dots,i_{k})=\emptyset$.
\endproclaim
\demo{Proof}
By definition
$$
I(i_0,\dots, i_{m+1}) = I(i_{m+1})\cap J(i_0,
\dots,i_{m}),
$$
so that
$$
I(i_0,\dots, i_{m+1}) \subseteq J(i_0,
\dots,i_{m}).
$$
Now $J(i_0,\dots,i_{m})$ is a closed set, so we have
$$
J(i_0,\dots, i_{m+1}) \subseteq \overline I(i_0,\dots, i_{m+1})
\subseteq J(i_0,\dots,i_{m}).
$$
Since this is true for all $m$, we easily obtain
$$
I(i_0,\dots,i_{k})\subseteq J(i_0,\dots,i_{l})
$$
for $0\le  l<k$.

To see that $J(i_0,\dots,i_{l})\subseteq
J(i_l)$, consider $t\in J(i_0,\dots,i_l)$. Every
neighbourhood of $t$ meets at least two members of the family
$\Cal I(i_0,\dots, i_l)$, and hence at least two members of the
family $\Cal I(i_l)$, so that $t\in J(i_l)$.

Finally, suppose that $i_m=i_l$ for some $0\le m<l\le k$. We have
$$
I(i_0,\dots, i_l)\subseteq I(i_l)\cap J(i_0,\dots ,i_m) \subseteq
I(i_l)\cap J(i_m),
$$
which is empty since $I(i)\cap J(i)=\emptyset $ for all $i$.
\enddemo

We shall be concerned especially with the sets $I(\bold i)$ when the
sequence $\bold i$ is strictly increasing. We shall write $\Sigma$
for the set of all such sequences $\bold i= (i_0, \dots, i_k)$ with
$k\ge 0$ and $i_0<i_1<\cdots <i_k$. We equip $\Sigma$ with a total
order $\prec$, defined by saying that $\bold i=(i_0, \dots,
i_k)\prec \bold j = (j_0, \dots, j_l)$ if {\it either} \roster
\item there exists $r\le \min\{k,l\}$ such that $i_s=j_s$ for
$0\le s<r$ and $i_r<j_r$, {\it or}
\item $k>l$ and $j_s=i_s$ for $0\le s\le l$.
\endroster
I am grateful to Gilles Godefroy who pointed out that this order
may be regarded as the usual lexicographic order if we think of
our finite sequences as infinite sequences terminating in a long
run of $\infty$'s.

\proclaim{Lemma 3.2}
Let $j= (j_0, \dots, j_l)\in \Sigma$ and write
$$
\align
A_1 = &\bigcup \Sb 0\le r\le l \\ j_{r-1}< i< j_r \endSb \overline I(j_0,\dots,
j_{r-1},i)\\
A_2 = &\bigcup \Sb k> l\\{i_k>i_{k-1}>\cdots>i_{l+1}>j_l}\endSb \overline I(j_0,\dots,
j_l,i_{l+1}, \dots, i_k)
\endalign
$$
Then
$$
\align
\bigcup_{\bold i\prec\bold j} \overline I(\bold i) &= A_1 \cup A_2\\
&= A_1 \cup J(\bold j).
\endalign
$$
In particular $\bigcup_{\bold i\prec\bold j} \overline I(\bold i)$
is a closed subset of $K$.
\endproclaim
\demo{Proof} It is clear that $A_2$ is exactly the union of the sets
$\overline I(\bold i)$ where $\bold i$ satisfies clause (2) in the
definition of the relation $\prec$.  If $\bold i=(i_0,\dots,i_k)$
satisfies clause (1) of that definition, then we have
$$
\overline I(i_0,\dots,i_k)\subseteq \overline I(i_0,\dots, i_r)=
\overline I(j_0,\dots,j_{r-1} i_r)\subseteq
A_1.
$$
It follows that $A_1$ is exactly the union of the sets $I(\bold
i)$ where $\bold i$ satisfies (2).

It is clear from the definitions that $A_2 \subseteq J(\bold j)$,
so, to prove the second equality, it will be enough to show that
$J(\bold j) \subseteq A_1 \cup A_2$.  Suppose then that $t\in
J(\bold j)$; for some $i$, we have $t\in I(i)$, and $ i$ is not
equal to any of the $j_s$,  since $J(\bold j)\subseteq J(j_s)$ and
$I(i)\cap J(i)=\emptyset.$
There are now two cases. If $i> j_l$ then
$$
t\in I(j_0,\dots,
j_l,i)\subseteq A_2.
$$
 If $i<j_l$ we choose $r$ minimal with
respect to $i< j_r$, noting that $i> j_{r-1}$, and observe that
$$
t\in I(j_0,\dots, j_{r-1}, i) \subseteq A_1.
$$
It is immediate that our set is closed, since we have shown it to
be the union of the closed set $J(\bold j)$ with finitely many
closures $\overline I(\bold i)$.
\enddemo

Given $\bold j$ and a finite subset $\Cal M$ of $\Cal I(\bold j)$
we shall write
$$
G(\bold j, \Cal M) = K \setminus\left( \bigcup_{\bold i\prec\bold j}
\overline I(\bold i)\cup \overline {\bigcup \Cal I(\bold j)\setminus
\Cal M}\right ),
$$
noting that this is an open subset of $K$.

Finally, we have a lemma which needs compactness of the space $K$.

\proclaim{Lemma 3.3}
Let $K$ be a compact space and let $\Cal S= \bigcup_{i\in
\omega}\Cal I(i)$ be a  covering of $K$ which is the union of regular
isolated families $\Cal I(i)$.
Let $H$ be a nonempty closed subset of $K$.  Then there exists a
minimal $\bold j\in \Sigma$ with $H\cap \overline I(\bold j) \ne
\emptyset$.  Moreover, $H\cap \overline I(\bold j)\subseteq I(\bold j)$
and there is a unique nonempty, finite $\Cal M\subseteq
\Cal I(\bold j)$ such that  $H\cap M\ne \emptyset$ for
all $M\in \Cal M$ and $H\subseteq G(\bold j, \Cal M)$.
\endproclaim
\demo{Proof}
If no minimal $\bold j$ exists, then we may find a strictly
decreasing sequence
$$
\bold j(0)\succ \bold j(1)\succ \cdots
$$
such that $H\cap \overline I(\bold j(n))\ne \emptyset $ for all $n$,
but such that $H\cap \overline I(\bold i)=\emptyset$ if $\bold
i\prec \bold j(m)$ for all $m$.  By Lemma 3.2, each of the sets
$H_n=H\cap \bigcup_{\bold i\prec \bold j(n)}\overline I(i)$ is
closed, and $H_n$ is nonempty since $H\cap \overline I(\bold
j_{n+1})\subseteq H_n$.  Hence, by compactness, the intersection
$\bigcap_{n\in\omega}H_n$ is nonempty.  But this means that $H\cap
\overline I(\bold i)\ne \emptyset$ for some $\bold i$ satisfying
$\bold i\prec \bold j(n)$ for all $n$.  This is a contradiction.

Working now with our minimal $\bold j$, we have $\overline I(\bold
j) = I(\bold j) \cup J(\bold j)$ and by Lemma~3.2 $J(\bold
j)\subseteq \bigcup_{\bold i\prec\bold j} \overline I(\bold i)$.
Thus, by minimality of $\bold j$, $H\cap J(\bold j)=\emptyset$ and
so $H\cap \overline I(\bold j) = H\cap I(\bold j)$.  The compact set
$H\cap \overline I(\bold j) $ is thus covered by the family $\Cal
I(\bold j)$, the elements of which are disjoint and open, relative
to $\overline I(\bold j)$.  Thus, if we define $\Cal M=\{M\in \Cal
I(\bold j): M\cap H\ne \emptyset\}$, it must be that $\Cal M$ is
finite. Finally, to see that $H\subseteq G(\bold j, \Cal M)$, we use
minimality of $\bold j$ again, and the observation that $\overline
{\bigcup \Cal I(\bold j)\setminus \Cal M}\subseteq \overline I(\bold
j)\setminus \bigcup \Cal M$, while $H\cap \overline I(\bold j)
\subseteq \bigcup \Cal M$.
\enddemo

When $K$ is a descriptive compact space having property $P$ with some metric $d$
then there exists, for each natural number $l$, a $\sigma$-isolated covering
$\Cal S^l=\bigcup_{i\in \omega}\Cal I^l(i)$ of $K$, consisting of
sets that are of $d$-diameter at most $2^{-l}$.  When $d$ is lower
semi-continuous, the sets $\widetilde N$ defined at the start of this
section are also of diameter at most $2^{-l}$.  In general, this
is not the case: however, each $\widetilde N$ is contained in the
$\Cal T$-closure of some set (namely $N$) of $d$-diameter at most
$2^{-l}$. We may summarize the situation in the form of a
proposition.

\proclaim{Proposition 3.4}
Let $(K,\Cal T)$ be a compact space equipped with a metric $d$
such that property $P(d, \Cal T)$ holds.  Then, for each $l\in
\omega$, there is a covering $\Cal S^l$ of $K$, which is the
union $\bigcup_{i\in \omega}\Cal I^l(i)$ of regular isolated families
$\Cal I^l(i)$, such that each $N\in \Cal S^l$ is contained in the
$\Cal T$-closure of some set of $d$-diameter at most $2^{-l}$
\endproclaim

From now on, we shall assume $S^l=\bigcup_{i\in \omega}I(i)$ to be as above, and shall
construct the associated $\Cal I^l(\bold i)$, $I^l(\bold i)$,
$J^l(\bold i)$ and $G^l(\bold i, \Cal M)$ as described in this
section.

\heading
4. Construction of a norm on $\Cal C(K)$
\endheading

We now set about constructing a norm on $\Cal C(K)$ when $K$ is a
descriptive compact space.  As well as the topological machinery set up in the
last section, we shall need one more ingredient.  Let $L$ be a
closed subset of $K$, let $l$ be a natural number, let $m,n$ be positive integers and let
$\bold i,\bold j\in \Sigma$; we write $\Cal B(L,l,m,n,\bold i,
\bold j)$ for the set of all pairs $(\Cal M, \Cal N)$ of finite
subsets of $\Cal I^l(\bold i), \Cal I^l(\bold j)$, respectively,
which satisfy $\# \Cal M=m$, $\# \Cal N=n$, $M\cap L\ne \emptyset$
for all $M\in \Cal M$, $N\cap L\ne \emptyset$
for all $N\in \Cal N$,  and
$$
\overline{\bigcup \Cal M} \cap \overline{\bigcup \Cal N} =
\emptyset.
$$

If $f\in \Cal C(K)$ and $L, \Cal M, \Cal N$ are as above, we set
$$
\Phi(f,L,\Cal M, \Cal N) = {\textstyle\frac12}\left(n^{-1} \sum_{N\in \Cal N}
\max f[L\cap\overline N]-m^{-1}\sum_{M\in \Cal M}
\min f[L\cap \overline M]\right ) ^+,
$$
noticing that $\Phi$ is a non-negative, positively homogeneous, convex function of its
argument $f$ and that
$$
\Phi(f,L,\Cal M, \Cal N) \le {\textstyle\frac12}\osc(f\restriction
L) \le \|f\restriction L\|_\infty.
$$

Whenever $(\Cal M, \Cal N)$ is a pair of finite sets as above, satisfying
$$
\overline{\bigcup \Cal M} \cap \overline{\bigcup \Cal N} =
\emptyset,
$$
we fix, once and for all, a
pair of closed subsets $(X(\Cal M, \Cal N), Y(\Cal M, \Cal N))$
such that
$$
X(\Cal M, \Cal N)\cup Y(\Cal M, \Cal N)=K, \quad X(\Cal M, \Cal
N)\cap \bigcup\Cal N= Y(\Cal M, \Cal
N)\cap \bigcup\Cal M=\emptyset.
$$

In the definition of our norm, we shall also need to fix positive
real numbers $c(\bold i)$ $(\bold i\in \Sigma)$ with $\sum_{\bold
i\in \Sigma}c(\bold i)\le 1$. We could, for instance, take
$$
c(i_0,i_1,\dots,i_k)= 2^{-2^{i_0}-2^{i_1}-\cdots -2^{i_k}}.
$$

\proclaim{Proposition 4.1}
There are unique non-negative real-valued functions $\Omega(f,L,l)
$,\newline
$\Theta(f,L,l,\bold i,\bold j,m,n) $,
$\Theta_p(f,L,l,\bold i,\bold j,m,n) $, $\Theta_p(f,L,l,\Cal M,
\Cal N)$ and
$\Psi(f,L, l, \Cal M, \Cal N)$, defined for functions $f\in
\Cal C(K)$, closed subsets $L$ of $K$, natural numbers $l,m,n,p$,
elements
$\bold i, \bold j$ of $ \Sigma$, and $(\Cal M, \Cal N)\in \Cal B(L,l,m,n,\bold i,
\bold j)$, which are convex  in their argument $f$, and which satisfy the
inequalities
$$
\align \Omega(f,L,l),\ &\Theta(f,L,l,\bold i, \bold j,m,n),
\Theta_p(f,L,l,\bold i, \bold j,m,n),\\\quad & \Theta_p(f,L,l,\Cal
M, \Cal N), \  \Psi(f,L, l,\Cal M, \Cal N)\le \|f\restriction
L\|_\infty,
\endalign
$$
as well as the relations
$$
\align
6\Omega(f,L,l)^2 &= \|f\restriction L\|_\infty^2 + \osc(f\restriction L)^2  \\
&\qquad\qquad +\sum_{\bold
i,\bold j\in \Sigma}c(\bold i)c(\bold
j)\sum_{m=1}^\infty\sum_{n=1}^\infty 2^{-m-n} \Theta(f,L,l,\bold i, \bold j,m,n)^2\\
\Theta(f,L,l,\bold i, \bold j,m,n)^2 &= \sum_{p=1}^\infty 2^{-p}\Theta_p(f,L,l,\bold i, \bold
j,m,n)^2\\
\Theta_p(f,L,l,\bold i, \bold j,m,n) &= \sup_{(\Cal M,\Cal N)\in
\Cal B(L,l,m,n,\bold i, \bold j)} \Theta_p(f,L,l,\Cal M, \Cal N)\\
2\Theta_p(f,L,l,\Cal M, \Cal N)^2 &= \Phi(f,L,l,\Cal M, \Cal N)^2 + p^{-1}\Psi(f,L,l,\Cal M, \Cal
N)^2\\
3\Psi(f,L,l,\Cal M, \Cal N)^2 &=\Omega(f,L\cap X(\Cal M,\Cal N),l)^2+
\Omega(f,L\cap Y(\Cal M,\Cal N),l)^2
\endalign
$$
We may define a norm $\|\cdot\|$ on $\Cal C(K)$, equivalent to
the supremum norm, by setting
$$
\|f\|^2 = \sum_{l=1}^\infty 2^{-l-1} \Omega(f,K,l)^2.
$$
\endproclaim
\demo{Proof} The functions $\Theta$ and $\Theta_p$ are defined in
terms of $\Psi$ and the known function $\Phi$ defined earlier. Hence
all we have to show is that the mutual recursion in the definitions
of $\Omega$ and $\Psi$ really does define something.  We do this by
applying a fixed-point theorem, as in [8].

Let $Z$ be the set of all tuples $(f,L,l,\Cal M,\Cal N)$ with
$f\in \Cal C(K)$, $L$ a closed subset of $K$, $l$ a positive
integer and $(\Cal M,\Cal N)\in \bigcup_{m,n,\bold i,\bold
j}B(L,l,m,n,\bold i,\bold j)$.  Let $\Cal Z$ be the set of all
pairs $(\Omega, \Psi)$ of non-negative real-valued functions $\Omega(f,L,l)$,
$\Psi(f,L,l,\Cal M,\Cal N)$, which are convex, symmetric and positively homogeneous
in their argument $f$, and which satisfy the inequalities
$$
\Omega(f,L,l),\ \Psi(f,L,l,\Cal M,\Cal N)\le \|f\|_\infty.
$$
Define a metric $\rho$ on $\Cal Z$ by setting
$$
\align
&\rho((\Omega,\Psi),(\Omega',\Psi'))= \\
\sup&\max\left\{|
\Omega(f,L,l)^2-\Omega'(f,L,l)^2|, |\Psi(f,L,l,\Cal M,\Cal N)^2-
\Psi'(f,L,l,\Cal M,\Cal N)^2|\right\},
\endalign
$$
where the supremum is taken over all $L,l,\Cal M,\Cal N$ and all
$f$ with $\|f\|_\infty\le 1$.  It is clear that this makes $\Cal
Z$ a complete metric space.

Now define a mapping $F:\Cal Z\to
\Cal Z$ by setting $F(\Omega,\Psi)=(\widetilde \Omega,\widetilde \Psi)$,
where
$$
3\widetilde\Psi(f,L,l,\Cal M, \Cal N)^2 =\Omega(f,L\cap X(\Cal M,\Cal N),l)^2+
\Omega(f,L\cap Y(\Cal M,\Cal N),l)^2,
$$
and
$$
\align
6\widetilde\Omega(f,L,l)^2 &= \|f\restriction L\|_\infty^2 + \osc(f\restriction L)^2  \\
&\qquad\qquad +\sum_{\bold
i,\bold j\in \Sigma}c(\bold i)c(\bold
j)\sum_{m=1}^\infty\sum_{n=1}^\infty 2^{-m-n} \Theta(f,L,l,\bold i, \bold
j,m,n)^2,
\endalign
$$
the function $\Theta$ being obtained from $\Psi$ via the formulae
in the statement of the Proposition.  It may be noted that, though
the function $\Theta$ is not symmetric in $f$, we do have
$$
\Theta(-f,L,l,m,n,\bold i,\bold j)=\Theta(f,L,l,n,m,\bold j,\bold
i),
$$
so that $\widetilde \Omega$ is symmetric.

It is easy to check that $\rho(F(\Omega,\Psi),F(\Omega',\Psi'))\le
\frac23\rho((\Omega,\Psi),(\Omega',\Psi'))$, so that $F$ has a
unique fixed point, by Banach's fixed point theorem.  This fixed
point yields the functions that we want, and hence enables us to
define the norm $\|\cdot\|$.
\enddemo

It is the norm defined in Proposition 4.1 that we shall
show to locally uniformly rotund in the case where $d$ is a
lower semi-continuous (or, more generally, Reznichenko)
metric fragmenting the descriptive compact space $K$.  By the
discussion at the end of \S 3 it will be enough to prove the
following theorem.

\proclaim{Theorem 4.2}
Let $(K,\Cal T)$ be a descriptive  compact space and let $d$ be a metric on
$K$ such that property $P(d,\Cal T)$ holds.  Let the norm $\|\cdot\|$ be defined as in
Proposition 4.1.  If $f$ be a function in $\Cal C(K)$ which is
$d$-uniformly continuous then the norm $\|\cdot\|$  locally uniformly convex at $f$.
\endproclaim

The proof of this theorem will occupy the remainder of the paper.
We shall consider a sequence $(f_r)$ in $\Cal C(K)$ which
satisfies
$$
{\textstyle\frac12}\|f\|^2+{\textstyle\frac12}\|f_r\|^2-\|{\textstyle\frac12}(f+f_r)\|^2\to
0,
$$
as $r\to \infty$. In the language introduced earlier, we are assuming that
the LUR hypothesis holds for $\|\cdot\|$ (and our given $f$ and
$f_r$).  We have to prove  that $f_r$ converges to
$f$ uniformly on $K$.  Given $\epsilon>0$, we may use uniform continuity of
$f$ to choose a positive integer $l$ such that
$$
d(t,u)\le 2^{-l}\implies |f(t)-f(u)|\le
{\textstyle\frac13}\epsilon.
$$

\proclaim{Lemma 4.3}
If $N\in \Cal S^l$ then the oscillation of $r$ on $N$ is at most
$\frac13\epsilon$.
\endproclaim
\demo{Proof}
As in Proposition 3.4, we are supposing that for each $N\in \Cal S^l$
there is some set $M$ of $d$-diameter at most $2^{-l}$ such that
$N$ is contained in the $\Cal T$-closure of $M$.  The uniform
continuity estimate tells us that the oscillation of $r$ on $M$ is
at most $\frac13 \epsilon$ and the $\Cal T$-continuity of $f$
enables us to extend this to $N$.
\enddemo

The definition of
our norm as an $\ell^2$-sum
$$
\|f\|^2 = \sum_{k=1}^\infty 2^{-k-1}\Omega(f,K,k)^2
$$
implies, thanks to an observation we made earlier, that the LUR
hypothesis holds for each of the functions $\Omega(\cdot,K,k)$ and
in particular for $\Omega(\cdot, K,l)$. This is all we shall use
in our proof that $\|f-f_r\|_\infty$ is eventually smaller than
$\epsilon$.

\heading 5. Good choices. \endheading
Let $L$ be a closed subset of of $K$, let $m,n$ be positive
integers and let $\bold i,\bold j\in \Sigma$. (Recall that $f$,
$\epsilon$ and $l$ are now fixed.)  For a pair
$(\widetilde{\Cal M},\widetilde{\Cal N})\in B(L,l,m,n,\bold i,\bold j)$, we define the
following non-negative real numbers:
$$
\align
A &= \min f[L]\\
a &= \max_{M\in \widetilde{\Cal M}} \inf f[L\cap \overline M]\\
\alpha &= \min f[L\setminus G^l(\bold i,\widetilde{\Cal M})]\\
\beta &= \max f[L\setminus G^l(\bold j,\widetilde{\Cal N})]\\
b &= \min_{N\in \widetilde{\Cal N}} \sup f[L\cap\overline N]\\
B &= \max f[L].
\endalign
$$
Of course, we have $a\ge A$, $\alpha\ge A$, $b\le B$ and $\beta\le
B$.  We shall say that the pair $(\widetilde{\Cal M},\widetilde{\Cal N})$ is a {\it good
choice} (of type $(m,n,\bold i,\bold j)$) on $L$ if
$$
\align
n^{-1}(B-\beta) &> (B-b)+(a-A)\quad\text{and}\\
m^{-1}(\alpha-A) &> (B-b)+(a-A).
\endalign
$$

\proclaim{Lemma 5.1}
If $L$ is a closed subset of $K$ and the oscillation of $f$ on $L$
is at least $\epsilon$ then there is at least one good choice on $L$.
\endproclaim
\demo{Proof}
Let $H_1=\{t\in L: f(t)=\max f[L]\}$ and apply Lemma 3.3.  There
exist $\bold j\in \Sigma$ and a finite subset $\widetilde{\Cal N}$ of $\Cal
I^l(\bold j)$ such that $H_1\cap N\ne \emptyset$ for all $N\in \widetilde{\Cal
N}$ and $H_1\subseteq G^l(\bold j,\widetilde{\Cal N})$.  It follows that, in
the notation just established, we have $B=b$ and $B>\beta$. A
similar argument applied to the set $H_2=\{t\in L: f(t)=\min
f[L]\}$ yields $\bold i$ and $\widetilde{\Cal M}$ such that $A=a$ and
$A<\alpha$.  To finish showing that $(\widetilde{\Cal M},\widetilde{\Cal N})$ is a good
choice, we need to check that $(\widetilde{\Cal M},\widetilde{\Cal N})$ is in $B(L,l,m,n,\bold i,\bold j)$,
and what remains to be proved is that $\overline M\cap \overline
N=\emptyset$ for all $M\in \widetilde{\Cal M}$ and all $N\in \widetilde{\Cal N}$.

Our choice of $l$ ensures that the oscillation of $f$ on each
$M\in \widetilde{\Cal M}$ and on each $N\in \widetilde{\Cal N}$ is at most $\epsilon/3$,
and, by continuity of $f$, the same holds for each $\overline M$
and each $\overline  N$.  Hence
$$
\align
\max f[\overline M] &\le A+\epsilon/3\\
\min f[\overline N] &\ge B-\epsilon/3,
\endalign
$$
for all such $M,N$.  Since we are assuming that the oscillation
$B-A$ of $f$ on $L$ is at least $\epsilon$, we deduce that
$\overline M\cap \overline N=\emptyset$ as claimed.
\enddemo

\proclaim{Lemma 5.2}
Let $L_1,L_2,\dots$ be a decreasing sequence of non-empty closed
subsets of $K$ with intersection $L$.  If $(\widetilde{\Cal M}, \widetilde{\Cal N})$ is a
good choice on $L$, then it is a good choice on $L_s$ for
arbitrarily large values of $s$.
\endproclaim
\demo{Proof}
Let us define $A,B,\alpha,\beta,b,B$ as above and set
$$
\align
A_s &= \min f[L_s]\\
a_s &= \max_{M\in \widetilde{\Cal M}} \inf f[L_s\cap\overline M]\\
\alpha_s &= \min f[L_s\setminus G^l(\bold i,\widetilde{\Cal M})],
\endalign
$$
with analogous definitions for $\beta_s,b_s,B_s$. Standard
compactness arguments show that $A_s\to A$ as $s\to \infty$, and
so on.  Hence the inequalities defining a good choice for $L_s$ do
hold for all sufficiently large $s$.
\enddemo

The third lemma in this section reveals why good choices are so
named: it is a ``rigidity condition'' of a type that occurs
commonly in LUR proofs.  It will be convenient to state it in
terms of ``strong attainment'' of a certain supremum, a notion
with which most readers will be familiar, but which we shall
nonetheless define explicitly.  If $(\gamma_i)_{i\in I}$ is a
bounded family of real numbers, we shall say that the supremum
$\sup_{i\in I}\gamma_i$ is {\it strongly attained} at $j$ if
$\sup_{i\in I\setminus\{j\}}\gamma_i<\gamma_j$.  This of course
implies that if $(i_r)$ is a sequence in $I$ and $\gamma_{i_r}\to
\sup_{i\in I}\gamma_i$ as $r\to \infty$, then $i_r=j$ for all
large enough $r$.

\proclaim{Lemma 5.3}
Let $L$ be a closed subset of $K$ and suppose that there exists a good choice
$(\widetilde{\Cal M},\widetilde{\Cal N})$ of type $(m,n,\bold i, \bold j)$ on $L$.  Then
the supremum $\sup\{\Phi(f,\Cal M,\Cal N):(\Cal M, \Cal N)\in
B(L,l,m,n,\bold i,\bold j)\}$ is strongly attained at $(\widetilde{\Cal M},\widetilde{\Cal N})$.
\endproclaim
\demo{Proof}
Let us write $A,a,\alpha,\beta,b,B$ for the quantities associated
with $(\widetilde {\Cal M},\widetilde{\Cal N})$ in the definition of a good
choice. So we have $\min f[L\cap \overline M] \le a$ and $\max
f[L\cap \overline N]\ge b$ for all $M\in \widetilde{\Cal M}$ and all
$N\in \widetilde{\Cal N}$.  Thus $2\Phi(f,\widetilde {\Cal M},\widetilde{\Cal N})\ge
b-a$.

Now suppose that $(\Cal M,\Cal N)$ is in $B(L,l,m,n,\bold i,\bold
j)$ and that $\Cal M\ne \widetilde{\Cal M}$.  Since $\Cal M$ and
$\widetilde{\Cal M}$ have the same number of elements, namely $m$, $\Cal
M$ must have at least one element $M_0$ which is not in $\widetilde{\Cal
M}$.  It follows from the definition of $G(\bold i,\widetilde{\Cal M})$
that $\overline M_0\cap G(\bold i,\widetilde{\Cal M})=\emptyset$ so that
$\min
f[L\cap \overline M_0]\ge \alpha$.  For the other $m-1$ elements
of $\Cal M$, we certainly have $\min
f[L\cap \overline M_0]\ge A$, and of course $\max f[L\cap
\overline N]\le B$ for all $N\in \Cal N$.  Hence
$$
2\Phi[f,\Cal M,\Cal N] \le \frac1n nB-\frac1m(\alpha+(m-1)A) =
B-A-\frac1m(\alpha-A).
$$
By the definition of a good choice, this is strictly smaller than
$b-a$, Similarly, we show that if $\Cal N\ne \widetilde{ \Cal N}$ then
$$
2\Phi[f,\Cal M,\Cal N] \le B-A-\frac1n(B-\beta),
$$
another quantity which is known to be smaller than $b-a$.
\enddemo

\heading
6. An application of Deville's Lemma
\endheading

We record for convenience the following version of Lemma VII.1.1 of
[2].

\proclaim{Lemma 6.1}
Let $(\phi_i)_{i\in I}$ and $(\psi_i)_{i\in I}$ be two pointwise-bounded
families of
non-negative, real-valued, convex functions on a real vector space
$Z$.    For $i\in I$ and positive integers $p$ define
functions $\theta_{i,p}$, $\theta_p$ and $\theta$ by setting
$$
\align
2\theta_{i,p}(x)^2 &= \phi_i(x)^2 + p^{-1}\psi_i(x)^2\\
\theta_p(x) &= \sup_{i\in I}\theta_{i,p}(x)\\
\theta(x)^2 &= \sum_{p=1}^\infty 2^{-p}\theta_p(x)^2.
\endalign
$$
Let $x$ and $x_r$ $(r\in \omega)$ be  elements of $Z$
and assume that
$$
{\textstyle\frac12}
\theta(x)^2+{\textstyle\frac12}\theta(x_r)^2-\theta({\textstyle\frac12}(x+x_r))^2\to
0
$$
as $r\to \infty$.  Then there is a sequence $(i_r)$ of elements of
$I$ such that
$$
\align
\phi_{i_r}(x)&\to \sup_{i\in I} \phi_i(x)\quad\text{and}\\
{\textstyle\frac12}\psi_{i_r}(x)^2 &+
{\textstyle\frac12}\psi_{i_r}(x_r)^2-\psi_{i_r}({\textstyle\frac12}(x+x_r))^2\to
0
\endalign
$$
as $r\to \infty$.
\endproclaim

\proclaim{Corollary 6.2}
If, in addition to the hypotheses of Lemma 6.1, we assume that the
supremum $\sup_{i\in I} \phi_i(x)$ is strongly attained at $j$,
then we may conclude that
$$
{\textstyle\frac12}\psi_{j}(x)^2 +
{\textstyle\frac12}\psi_{j}(x_r)^2-\psi_{j}({\textstyle\frac12}(x+x_r))^2\to
0.
$$
\endproclaim
\demo{Proof}
This is of course automatic, since the assumptions imply that a
sequence $(i_r)$ for which
$$
\phi_{i_r}(x)\to \sup_{i\in I} \phi_i(x)
$$
as $r\to\infty$ must necessarily satisfy $i_r=j$ for all large
enough $j$.
\enddemo

We may rephrase the statement of this corollary by saying that if the LUR
hypothesis holds for $\theta$ and the supremum
$\sup_{i\in I} \phi_i(x)$ is strongly attained at $j$, then the
LUR hypothesis holds for $\psi_j$.  It is precisely this
formulation that we shall be applying in the next result, where we
return to the proof of Theorem 4.2 and where of course we are
still dealing with fixed $f,f_r,\epsilon$ and $l$.

\proclaim{Proposition 6.3}
Let $L$ be a closed subset of $K$ and assume that the LUR
hypothesis holds for $\Omega(\cdot,L,l)$.  If $(\widetilde{\Cal
M},\widetilde{\Cal N})$ is a good choice on $L$ then the LUR hypothesis
holds for $\Omega(\cdot,L\cap X(\widetilde{\Cal M},\widetilde{\Cal N}),l)$
and $\Omega(\cdot,L\cap Y(\widetilde{\Cal M},\widetilde{\Cal N}),l)$.
\endproclaim
\demo{Proof}
Let $(\widetilde{\Cal M},\widetilde{\Cal N})$ be of type $(m,n,\bold i,\bold
j)$.  The expression for $\Omega(\cdot,L,l)$ as an $\ell^2$-sum
implies that the LUR hypothesis holds for
$\Theta(\cdot,L,l,m,n,\bold i,\bold j)$, which is readily
recognizable as a function to which we may apply Deville's lemma.
Moreover, by Lemma 5.3, we are in the situation where the supremum
$\sup_{(\Cal M,\Cal N)\in B(m,n,\bold i,\bold j)}\Phi(f,L,l,\Cal
M,\Cal N)$ is strongly attained at $(\widetilde{\Cal M},\widetilde{\Cal N})$.
So, by the above corollary, the LUR hypothesis holds for
$\Psi(\cdot,L,l,\widetilde{\Cal M},\widetilde{\Cal N})$.  The formula for this
as an $\ell^2$ sum now shows that the LUR hypothesis holds for
$\Omega(\cdot,L\cap X(\widetilde{\Cal M},\widetilde{\Cal N}),l)$
and $\Omega(\cdot,L\cap Y(\widetilde{\Cal M},\widetilde{\Cal N}),l)$, as
claimed.
\enddemo

\heading
7. Putting together the pieces
\endheading

\proclaim{Lemma 7.1} Let $L$ be a closed subset of $K$ on which the
oscillation of $f$ is smaller than $\epsilon$.  If the LUR
hypothesis holds for $\Omega(\cdot,L,l)$ then $\|(f-f_r)\restriction
L\|_\infty<\epsilon$ for all large enough $r$.
\endproclaim
\demo{Proof} From the formula for $\Omega$ as an $\ell^2$-sum, we
see that the LUR hypothesis holds for the the convex functions
$g\mapsto \|g\restriction L\|_\infty$ and $g\mapsto
\osc(g\restriction L)$. So in particular, $\|f_r\restriction
L\|_\infty\to\|f\restriction L\|_\infty$,
$\|{\textstyle\frac12}(f+f_r)\restriction
L\|_\infty\to\|f\restriction L\|_\infty$ and $\osc(f_r\restriction
L)\to\osc(f\restriction L)$ as $r\to \infty$. The required result
follows from a fairly standard argument.
 Let us write $\osc(f\restriction
L)=\epsilon-4\eta$ and suppose that $r$ is large enough for us to
have
$$
\align
\|f_r\restriction L\|_\infty &< \|f\restriction L\|_\infty +\eta\\
\|{\textstyle\frac12}(f+f_r)\restriction L\|_\infty &>
\|f\restriction L\|_\infty
-\eta\\
\osc(f_r\restriction L) &<\epsilon - 3\eta.
\endalign
$$
There exists $t\in K$ with
$|{\textstyle\frac12}(f+f_r)(t)|>\|f\restriction L\|_\infty -\eta$,
and we may assume that $(f+f_r)(t)>0$.  It follows that
$$
\align f_r(t) &> 2\|f\restriction L\|_\infty
-2\eta -\|f\restriction L\|_\infty\\
&= \|f\restriction L\|_\infty-2\eta\\
f(t) &> 2\|f\restriction L\|_\infty
-2\eta -\|f_r\restriction L\|_\infty\\
&> \|f\restriction L\|_\infty -3\eta
\endalign
$$
Now for any $u\in L$ we have
$$
\align
f(u) &\ge f(t)-\osc(f\restriction L)\\
     &> \|f\restriction L\|_\infty
-3\eta -\epsilon + 4\eta\\
  &= \|f\restriction L\|_\infty -\epsilon + \eta\\
f(u) &\le\|f\restriction L\|_\infty\\
f_r(u) &\ge f_r(t)-\osc(f_r\restriction L)\\
   &> \|f\restriction L\|_\infty-2\eta-\epsilon + 3\eta\\
    &= \|f\restriction L\|_\infty-\epsilon + \eta\\
f_r(u) &\le \|f_r\restriction L\|_\infty\\
&<\|f\restriction L\|_\infty+\eta.\\
\endalign
$$
It follows immediately that $|f(u)-f_r(u)|<\epsilon$.
\enddemo

\proclaim{Proposition 7.2}
There is a finite covering $\Cal L$ of $K$ with closed subsets
such that the LUR hypothesis holds for $\Omega(\cdot,L,l)$
and the oscillation of $f$ on $L$ is smaller than $\epsilon$,
for each $L\in \Cal L$.
\endproclaim
\demo{Proof}
We shall define a tree $\Upsilon$ whose elements will be certain
pairs $(L,s)$ with $L$ a closed subset of $K$ and $s$ a natural
number. We shall give a recursive definition which will specify
which such pairs are nodes of our tree, and shall define the tree
ordering by saying which (if any) nodes are the {\it immediate
successors} of a given $(L,s)$.  To do this, we shall need
to fix a mapping $\tau:\omega\to
\omega\times\omega\times\Sigma\times \Sigma$ with the property
that each quadruple $(m,n,\bold i,\bold j)$ occurs as $\tau(s)$
for infinitely many $s\in \omega$.

It will be ensured during the construction
that, whenever $(L,s)\in \Upsilon$, the LUR hypothesis holds for
$\Omega(\cdot,L,l)$.  We start  by declaring that
there is one minimal node $(K,0)$. (Notice that our hypotheses do
ensure that the LUR hypothesis holds for $\Omega(\cdot,K,l)$.)  If
$(L,s)$ is a node of our tree then there are three possibilities:
\roster
\item if the oscillation of $f$ on $L$ is smaller than $\epsilon$
then $(L,s)$ has no immediate successors in the tree (that is to
say, $(L,s)$ is a maximal element);
\item if the oscillation of $f$ on $L$ is at least $\epsilon$ and
there is a good choice $(\Cal M,\Cal N)$ of type $\tau(s)$ on $L$
then we introduce into $\Upsilon$ two immediate successors, $(L\cap
X(\Cal M,\Cal N),s+1)$ and $(L\cap
X(\Cal M,\Cal N),s+1)$, of $(L,s)$  (Notice that, by Proposition 6.3,
the LUR hypothesis holds for the $\Omega$ functions
associated with these two new nodes.);
\item if the oscillation of $f$ on $L$ is at least $\epsilon$ but
no good choice of type $\tau(s)$ exists, then we introduce just
one immediate successor $(L,s+1)$ of $(L,s)$ into the tree.
\endroster

We shall now show that the tree $\Upsilon$ we have just constructed
has only finitely many elements.  By K\"onig's Lemma, it will be
enough to show that $\Upsilon$ has no infinite branch. So suppose,
if possible, that there is a sequence $(L_s)_{s\in \omega}$ of
closed subsets of $K$ such that the pairs $(L_s,s)$ are nodes of
$\Upsilon$ and such that, for each $s$, $(L_{s+1},s+1)$ is an
immediate successor of $(L_s,s)$ in $\Upsilon$. The sets $L_s$ form
a decreasing sequence of closed subsets of $K$; let us write $L$ for
their intersection. By a compactness argument, the oscillation
$\osc(f\restriction L_s)$ tends to $\osc(f\restriction L)$ as $s\to
\infty$.  Since each $(L_s,s)$ has successors in $\Upsilon$, we have
$\osc(f\restriction L_s)\ge \epsilon$ for each $s$, and we can thus
deduce that $\osc (f\restriction L)\ge \epsilon$. So, by Lemma 5.1,
there is a good choice $(\Cal M,\Cal N)$ on $L$, of type $(m,n,\bold
i,\bold j)$ say. By Lemma 5.2, $(\Cal M,\Cal N)$ is also a good
choice on $L_s$ for all sufficiently large $s$. Recalling that
$\tau(s)=(m,n,\bold i,\bold j)$ for infinitely many values of $s$,
we see that we can choose $s$ such that $(\Cal M,\Cal N)$ is a good
choice on $L_s$ of type $\tau(s)$.  The way we constructed the tree
$\Upsilon$ means that $L_{s+1}$ is one or other of the two sets
$L_s\cap X(\Cal M,\Cal N)$ and $L_s\cap Y(\Cal M,\Cal N)$. So one or
other of $L_{s+1} \cap \bigcup \Cal N$ and $L_{s+1}\cap \bigcup\Cal
M$ is empty.  But this is absurd, since $L_{s+1}\supseteq L$ and the
sets $L\cap M$,  $L\cap N$ are nonempty for all $M\in \Cal M$ and
$N\in \Cal N$.

Having proved that $\Upsilon$ is finite, we define
$\Cal L $ to be the set of all $L$ such that there is a maximal
element of $\Upsilon$ of the form $(L,r)$.  Our construction
ensures that the LUR hypothesis holds for $\Omega(\cdot,L,l)$ for
each such $L$, and, by maximality, the oscillation of $f$ on any
such $L$ is smaller than $\epsilon$.  We just need to show that
$\bigcup \Cal L=K$.  This is most easily proved by induction: for
each $s$ let $\Cal L_s=\{L: (L,s)\in \Upsilon\}$; I claim that,
for all $s$, $\bigcup \Cal L\cup \bigcup \Cal L_s=K$.  Certainly
this is true for $s=0$ since $\Cal L_0=\{K\}$.  To deal with the
inductive step let $t\in \bigcup \Cal L\cup \bigcup L_s$ be given.
If $t\in \bigcup\Cal L$ there is no problem, so assume that $t\in L$ for some $L\in \Cal
L_s$.  By the construction of $\Upsilon$, one of (1), (2) and (3)
occurs for the pair $(L,s)$.  If it is (1) then $t\in L\in \Cal
L$.  If it is (2) then $t$ is one or other of the two sets $L\cap
X$ and $L\cap Y$ which themselves are members of $\Cal L_{s+1}$.
Finally if it is (3) then $t\in L\in \Cal L_{s+1}$.  In all cases,
we have $t\in\bigcup \Cal L\cup \bigcup \Cal L_{s+1}$, which
completes our proof by induction.  Since $\Upsilon$ is finite,
$\Cal L_s$ is empty for large enough $s$, which shows that
$\bigcup\Cal L=K$.
\enddemo

We can now finish the proof of Theorem 4.2.  Indeed, by Proposition
7.2, $K$ is the union of finitely many subsets $L$, for each of
which $\sup_{t\in L}|f_r(t)-f(t)|$ is eventually smaller than $\epsilon$.
So $\|f_r-f\|_\infty$ is eventually smaller than $\epsilon$,
which is what we wanted to prove.

\enddocument